\documentclass[12pt]{article}
\usepackage{cite,amssymb,amstext,amsmath}

\def\author#1{{\it#1},}
\let\Bbb\mathbb
\let\boldsymbol\relax
\def\book#1{{\sl#1}}
\def\paper#1{#1,}
\let\Cal\mathcal

\long\def\quote#1{%
  {\bigskip\small
  \advance\leftskip 2.2\parindent
  \everydisplay{\hskip\leftskip}%
  \parindent 0pt
  \noindent\llap{\cite{S4}:\ \ }#1\par\bigskip}%
}

\newbox\GrossOneBox
\newbox\grossOneBox
\newbox\GrossOneBoxN
\newbox\grossOneBoxN

\setbox\GrossOneBox\hbox{\raise-.4pt\hbox{${\rm O}\mskip-10.2mu\raise-2.7pt\hbox{$^1$}\mskip3mu$}}
\setbox\grossOneBox\hbox{\raise-4.5pt\hbox{\small$\mskip-1.3mu{\hbox{\rm\footnotesize O}}\mskip-10.0mu\raise-3.0pt\hbox{$^1$}\mskip1mu$}}
\setbox\GrossOneBoxN\hbox{\raise-.2pt\hbox{\small${\rm O}\mskip-10.5mu\raise-2.7pt\hbox{$^1$}\mskip3mu$}}
\setbox\grossOneBoxN\hbox{\raise-4.5pt\hbox{\small$\mskip-1.3mu{\hbox{\rm\footnotesize O}}\mskip-10.0mu\raise-3.0pt\hbox{$^1$}\mskip1mu$}}

\def\GrossOne{{\copy\GrossOneBox}}
\def\grossOne{{\copy\grossOneBox}}
\def\GrossOneN{{\copy\GrossOneBoxN}}
\def\grossOneN{{\copy\grossOneBoxN}}

\begin{document}

\thispagestyle{empty}

\hbox{}
\vskip20mm
\begin{center}
  \LARGE
  A.\,E.\,Gutman\footnote{The work of the first author is supported by the Russian Science Support Foundation.}\\
  S.\,S.\,Kutateladze
  \vskip10mm
  {\bf A TRIVIAL FORMALIZATION\\
  OF THE THEORY OF GROSSONE}
\end{center}
\vskip20mm

{\leftskip5mm\rightskip5mm
\parindent0pt
{\bf Abstract. }
A~trivial formalization is given for the~informal reasonings presented in a~series of papers by Ya.\,D.\,Sergeyev
on a~positional numeral system with an~infinitely large base, grossone;
the~system which is groundlessly opposed by its originator to the~classical nonstandard analysis.

\bigskip


\bigskip

\par}

\vfill
{\parindent 110mm
\copyright\ A.\,E.~Gutman,~ 2008\par
\copyright\ \,S.\,S.~Kutateladze,~ 2008
\bigskip
}

\newpage

In~recent years Ya.\,D.\,Sergeyev has published a~series of papers~\cite{S1,S2,S3,S4,S5}
in which a~positional numeral system is advanced related to the~notion of grossone%
\footnote{
  \,The~term ``grossone'' belongs to Ya.\,D.\,Sergeyev, has no relevance to the~usual meaning of the~noun
  ``gross'' in~English, and stems most likely from ``gro\ss'' in~German or ``grosso'' in~Italian.
}.
Ya.\,D.\,Sergeyev opposes his system to nonstandard analysis
and regards the~former as resting on different mathematical, philosophical, etc.~doctrines.
The~aim of the~present note is to properly position the~papers by~Ya.\,D.\,Sergeyev
on developing numeral systems.
It~turns out that a~model of Ya.\,D.\,Sergeyev's system is provided by the~initial segment $\{1,2,\dots,\nu!\}$
of~the~nonstandard natural scale up~to the~factorial~$\nu!$ of an~arbitrary actual infinitely large natural~$\nu$.
Such a~factorial serves as a~model of Ya.\,D.\,Sergeyev's grossone, thus demonstrating the~place occupying
by the~numeral system he proposed.

As~the main source we have chosen~\cite{S4}, the~latest available paper by Ya.\,D.\,Sergeyev,
which contains a~detailed description of his basic ideas.

\quote{.~.~.~the approach used in this paper is different also with respect
to the~nonstandard analysis~.~.~. and built using Cantor's ideas.}

\noindent
In~the present note we are about to show that, contrary to what is expected by the~author of~\cite{S4},
his indistinct definitions of grossone and the~concomitant notions
admit an~extremely accurate and trivial formalization within the~classical nonstandard analysis.

\quote{The~infinite radix of the~new system is introduced as the~number of elements of the~set~${\Bbb N}$
of~natural numbers expressed by the~numeral ${\GrossOneN}$ called {\it grossone}.}

\noindent
Use the~formalism of the~internal set theory IST by E.\,Nelson~\cite{N}
or any of the~classical external set theories,
for instance, EXT by~K.\,Hrba\v{c}ek~\cite{H} or NST by~T.\,Kawai~\cite{K}
(see also the~monographs~\cite{KR,GKK}).
As~usual, ${}^\circ\!X$ denotes the~standard core of a~set~$X$,
i.e., the~totality of all standard elements of~$X$.
In~particular, ${}^\circ{\Bbb N}$ is the~totality of all finite (standard) naturals.
Fix~an arbitrary infinitely large natural~$\nu$ and denote its factorial by~${\GrossOne}$:
$$
{\GrossOne}=\nu!\,,\quad\text{where }\,\nu\in{\Bbb N},\ \nu\approx\infty.
$$
Show that ${\GrossOne}$ possesses all properties of ``grossone''
(postulated as well as implicitly presumed in~\cite{S4}).

A~possible approach to an~adequate formalization (in the~sense of~\cite{S4}) of the~notion of size or ``the~number of elements''
of an~arbitrary set~$A$ of standard naturals
(i.e., of an~external subset~$A\subset{}^\circ{\Bbb N}$)
consists in assigning the~natural $\|A\| = |{}^*\!A\cap\{1,2,\dots,{\GrossOne}\}|$ to each $A$,
where ${}^*\!A$ is the~standardization of~$A$ and $|X|$ is the~size
(in the~usual sense) of a~finite internal set~$X$.
In~this case it is clear that $\|{}^\circ{\Bbb N}\|={\GrossOne}$,
which agrees with the~fore-quoted ``definition'' of grossone.
Note also that, due to the~external induction, the~function $A\mapsto\|A\|$
possesses the~additivity property (presumed in~\cite{S4}):
$\bigl\|\bigcup_{k=1}^n A_k\bigr\| = \sum_{k=1}^n\|A_k\|$
for every family of pairwise disjoint sets $A_1,\dots,A_n\subset{}^\circ{\Bbb N}$, $n\in{}^\circ{\Bbb N}$.

Another approach (which is more trivial and considerably closer to that of~\cite{S4})
to defining the~number of elements consists in ``replacing'' the~set~${}^\circ{\Bbb N}$
with the~initial segment
$$
{\Cal N}=\{1,2,\dots,{\GrossOne}\}
$$
of the~natural scale
and considering the~usual size $|A|\in{\Bbb N}$
of each internal set $A\subset{\Cal N}$.
In~this case, again, $|{\Cal N}|={\GrossOne}$; and the~additivity of the~counting measure
$A\mapsto|A|$ needs no argument.

\newpage
\quote{The~new numeral ${\GrossOneN}$ allows us to write down the~set, ${\Bbb N}$, of natural numbers in the~form
$$
{\Bbb N} = \{1,\ 2,\ 3,\ \dots,\ {\GrossOneN}-2,\ {\GrossOneN}-1,\ {\GrossOneN}\,\}
$$
because {\it grossone has been introduced as the~number of elements of the~set of natural numbers}
(similarly, the~number~$3$ is the~number of elements of the~set $\{1,2,3\}$). Thus, grossone is the~biggest natural number~.~.~.}

\noindent
While crediting the~author of~\cite{S4} for the~audacious extrapolation of the~properties of the~number~$3$,
we~nevertheless cannot accept the~fore-quoted agreement if for no other reason than the~fact that
the~set~${\Bbb N}$ of naturals (in the~popular sense of this fundamental notion)
has no greatest element (with respect to the~classical order).
In~order to keep the~traditional sense for the~symbol~${\Bbb N}$
(and being governed by ``Postulate~3. {\it The~part is less than the~whole}'' of~\cite{S4}),
instead of reusing this symbol for the~proper subset $\{1,2,\dots,{\GrossOne}\}\subset {\Bbb N}$
we decided to give the~latter a~less radical notation,~${\Cal N}$.

\quote{{\it The\,Infinite\,Unit\,Axiom\/} consists of the~following three statements:
\par\smallskip
{\it Infinity}. For~any finite natural number $n$ it follows $n < {\GrossOneN}$.
\par\smallskip
{\it Identity}. The~following relations link ${\GrossOneN}$ to identity elements $0$ and $1$
$$\textstyle
0\mathbin{\boldsymbol\cdot}{\GrossOneN}={\GrossOneN}\mathbin{\boldsymbol\cdot}0=0,\quad
{\GrossOneN}-{\GrossOneN}=0,\quad
\frac{\grossOneN}{\grossOneN}=1,\quad
{\GrossOneN}^0=1,\quad
1^{\grossOneN}=1,\quad
0^{\grossOneN}=0.
$$
{\it Divisibility}. For~any finite natural number $n$ sets ${\Bbb N}_{k,n}$, $1\leqslant k\leqslant n$, being the~$n$th
parts of the~set, ${\Bbb N}$, of natural numbers have the~same number of elements
indicated by the~numeral~$\frac{\grossOneN} n$, where
$$
{\Bbb N}_{k,n} = \{k,\ k+n,\ k+2n,\ k+3n,\ \dots\},\quad
1\leqslant k\leqslant n,\quad
\bigcup_{k=1}^n {\Bbb N}_{k,n}={\Bbb N}.
$$}

\noindent
Since ${\GrossOne}=\nu!$ is an~infinitely large number, it satisfies {\it Infinity}.
Every natural meets {\it Identity\/}, and so does~${\GrossOne}$.
Presenting the~factorial of an~infinitely large number, ${\GrossOne}$ is divisible by every standard natural.
Moreover, if $n\in{}^\circ{\Bbb N}$, $1\leqslant k\leqslant n$, and
\begin{align*}
{}^\circ{\Bbb N}_{k,n}&=\{k+(m-1)n\,:\,m\in{}^\circ{\Bbb N}\},\\
{\Cal N}_{k,n}&={\Cal N}\cap\{k+(m-1)n\,:\,m\in{\Bbb N}\};
\end{align*}
then $\|{}^\circ{\Bbb N}_{k,n}\|=|{\Cal N}_{k,n}|=\frac {\grossOne} n$.
Hence, ${\GrossOne}$ meets {\it Divisibility}.

\quote{It~is worthy to emphasize that, since the~numbers $\frac{\grossOneN} n$
have been introduced as numbers of elements of sets ${\Bbb N}_{k,n}$, they are integer.}

\noindent
If~a~number is declared natural, it naturally cannot occur unnatural.
To~remove all doubts, we~suggest a~rigorous and detailed justification for satisfiability of the~above postulate:
for~every $n\in{}^\circ{\Bbb N}$ we have $n<\nu$ and thus
$$
\text{the number }\ \frac{\GrossOne} n=\frac{\nu!}n=
\frac{1\mathbin{\boldsymbol\cdot} 2\mathbin{\boldsymbol\cdot}\ldots
\mathbin{\boldsymbol\cdot} n\mathbin{\boldsymbol\cdot}\ldots\mathbin{\boldsymbol\cdot}\nu}n\ \text{ is integer}.
$$

\newpage
\quote{The~introduction of grossone allows us to obtain the~following interesting result:
the~set~${\Bbb N}$ is~not a~monoid under addition. In~fact, the~operation ${\GrossOneN}+1$ gives us
as the~result a~number grater than ${\GrossOneN}$. Thus, by definition of~grossone, ${\GrossOneN}+1$ does
not belong to ${\Bbb N}$ and, therefore, ${\Bbb N}$~is not closed under addition and is not a~monoid.}

\noindent
Indeed, ${\GrossOne}\in\{1,2,\dots,{\GrossOne}\}={\Cal N}$, but ${\GrossOne}+1\notin\{1,2,\dots,{\GrossOne}\}={\Cal N}$.
(However, taking it into account that ${\Cal N}$ is not the~set of all naturals,
the~above trivial observation is unlikely ``interesting.'')

\quote{.~.~.~adding the~Infinite\,Unit\,Axiom to the~axioms of natural numbers defines
the~set of extended natural numbers indicated as $\widehat{\Bbb N}$ and
including ${\Bbb N}$ as a~proper subset
$$
\widehat{\Bbb N} = \{1,\ 2,\ \dots,\ {\GrossOneN}-1,\ {\GrossOneN},\ {\GrossOneN}+1,\
\dots,\ {\GrossOneN}^2-1,\ {\GrossOneN}^2,\ {\GrossOneN}^2+1,\ \dots\}.
$$}
\noindent
In~fact, ${}^\circ{\Bbb N}$ and ${\Cal N}$ are both proper subsets of the~set~${\Bbb N}$ of {\it all\/} naturals.
(As is known, the~radical formalism of IST saves us from considering ``extended numbers.'')

\bigskip
We~permit ourselves to pass over other numerous descriptions of the~properties of grossone and the~accompanying notions in~\cite{S4},
since the~corresponding analysis is quite analogous to that above (and equally trivial).
However, we cannot help commenting the~declared elimination of~Hilbert's paradox of the~Grand Hotel:

\quote{.~.~.~it is well known that Cantor's approach leads to some ``paradoxes''~.~.~.
Hilbert's Grand Hotel has an~infinite number of rooms~.~.~.
If~a new guest arrives at the~Hotel where every room is occupied,
it is, nevertheless, possible to find a~room for him.
To~do so, it is necessary to move the~guest occupying room 1 to room 2, the~guest
occupying room 2 to~room 3, etc.
In~such a~way room 1 will be available for the~newcomer~.~.~.
\par\medskip
.~.~.~In the~paradox, the~number of the~rooms in the~Hotel is countable.
In~our terminology this means that it has ${\GrossOneN}$ rooms~.~.~.
Under the~Infinite Unit Axiom this procedure is not possible
because the~guest from room ${\GrossOneN}$ should be moved to room ${\GrossOneN}+1$ and the~Hotel has only ${\GrossOneN}$~rooms.
Thus, when the~Hotel is full, no more new guests can be accommodated~--- the~result
corresponding perfectly to Postulate 3 and the~situation taking place in~normal
hotels with a~finite number of rooms.}

\noindent
The~following unpretentious ``paradox of the~Gross Hotel'' is brought to the~audience's attention:
Even though all the~grossrooms $1,2,\dots,{\GrossOne}$ are occupied, it is easy
to accommodate one more client in the~Gross Hotel. To~this end it suffices
to move the~guest occupying room~$n$ to~room~$n+\nobreak1$ for each {\it finite\/}~$n$.
Since $n+1<{\GrossOne}$ for all finite~$n$,
all the~former guests get their rooms in the~Gross Hotel, while room 1 becomes free for a~newcomer.

\bigskip
Besides a~babbling theorization around grossone, \cite{S4} includes an~``applied'' part
dedicated to a~new positional numeral system with base~${\GrossOne}$.
(The system is meant for becoming a~foundation for ``Infinity Computer''~\cite{S5}
which is able to operate infinitely large and infinitesimal numbers.)
Unfortunately, the~corresponding exposition remains highly informal,
and even crucial definitions are substituted with allusions and illustrating examples.

\newpage
\quote{In~order to construct a~number $C$ in the~new numeral positional system with base ${\GrossOneN}$ we~subdivide~$C$ into groups
corresponding to powers of ${\GrossOneN}$:
$$
C = c_{p_m}{\GrossOneN}^{p_m} + \cdots + c_{p_1}{\GrossOneN}^{p_1} + c_{p_0}{\GrossOneN}^{p_0} +
c_{p_{-1}}{\GrossOneN}^{p_{-1}} + \cdots + c_{p_{-k}}{\GrossOneN}^{p_{-k}}.
$$
.~.~.~Finite numbers $c_i$ are called {\it infinite grossdigits\/}
and can be both positive and negative;
numbers $p_i$ are called {\it grosspowers} and can be finite, infinite, and infinitesimal
(the~introduction of infinitesimal numbers will be given soon).
The~numbers $p_i$ are such that $p_i > 0$, $p_0 = 0$, $p_{-i} < 0$ and
$$
p_m > p_{m-1} > \cdots > p_2 > p_1 > p_{-1} > p_{-2} > \cdots > p_{-(k-1)} > p_{-k}.
$$
.~.~.~Finite numbers in this new numeral system are represented by numerals having
only one grosspower equal to zero~.~.~.
\par\medskip
.~.~.~all grossdigits $c_i$, $-k\leqslant i\leqslant m$, can be integer or fractional~.~.~.
Infinite numbers in this numeral system are expressed by numerals having at~least one grosspower grater than zero~.~.~.
Numerals having only negative grosspowers represent infinitesimal numbers.}

\noindent
In~the fore-quoted definitions, combinations of the~terms ``finite,'' ``infinite,'' and ``number'' seem to be used quite vaguely.
For~instance, it is unclear from the~text whether a~numeral is assumed infinite (and in what sense) if it is not finite (in some sense).
Following the~definitions of~\cite{S4} literally, a~grosspower can be finite, infinite, and (or?)~infinitesimal,
while ``finite'' means~$c{\mskip1mu}{\GrossOne}^0$ (a~grossdigit~$c$, a~rational numeral),
``infinite'' is expressed by a~numeral having al least one strictly positive grosspower,
and ``infinitesimal'' is a~numeral whose grosspowers are all strictly negative.
Seemingly, this implies that a~grosspower cannot be equal to, say, ${\GrossOne}^0+{\GrossOne}^{-1}$,
but the~subsequent examples of~\cite{S4} show that this is not so, and arbitrary numerals can serve as~grosspowers.
In~addition, the~reason is completely unclear for choosing the~terms ``infinite'' and ``infinitesimal''
exactly for the~classes of numerals mentioned in the~quote.
For~instance, the~numeral~$a={\GrossOne}^{{\grossOne}{\raise-6pt\hbox{$^{^{-1}}$}}}$ \!(with grosspower ${\GrossOne}^{-1}>0$)
is ``infinite'' by definition, while, obviously, $1<a<2$.
On~the other hand, the~numeral $b={\GrossOne}^{{\grossOne}{\raise-6pt\hbox{$^{^{-1}}$}}}{\mskip-8mu}-1$
is also considered ``infinite'' and not ``infinitesimal,''
while, as is easily seen, $b$ is infinitely close to zero in the~sense that $-c<b<c$ for every finite~$c>0$.

Regardless of terminological discipline, the~fore-quoted definition of numerals~$C$
cannot be considered formal if for no other reason than the~participating notion of (``infinite'' and ``infinitesimal'')
grosspowers depends on the~initial notion of numeral, thus leading to a~vicious circle.
In~addition, from the~illustrations of~\cite{S4} it is clear that the~positional system proposed admit
syntactically different numerals with coincident values: for instance,
$0{\mskip1mu}{\GrossOne}^0\equiv0{\mskip1mu}{\GrossOne}^1$,
$1{\GrossOne}^0\equiv1{\GrossOne}^{\hbox{\footnotesize0}\mskip2mu{\raise4.3pt\hbox{\grossOneN}}^{\mskip1mu0}}$.
(The notion of the~value of a~term and the~equivalence relation~$\equiv$ are clarified in~\cite{GK}.)
At~the same time, \cite{S4} misses not only the~corresponding stipulations (easy to guess though)
but also any attempts of justifying the~unambiguity of the~positional system, even under implicit stipulations.
Observe also that the~description of~\cite{S4} for the~algorithms of calculating the~sum and product of numerals
(i.e.,~of finding the~corresponding equivalent numeral) is~very superficial, since it does not touch upon the~problem
of recognizing equivalent numerals (which is necessary for collecting similar terms)
and that of comparing them (which is necessary for collating the~summands in order of their ``grosspowers'').
It~is thus not surprising that the~patent application~\cite{S5} reports on the~development of ``Infinity Calculator''
which is able to handle numerals admitting ``finite exponents'' only.

To~provide some justification, we briefly described in~\cite{GK}
one of the~possible approaches to~formalization of the~notion of numeral
as well as the~corresponding algorithmic procedures.

\vfill

\noindent
{\it Alexander E.~Gutman}
\par\smallskip
{\leftskip\parindent\small
\noindent
Sobolev Institute of Mathematics\\
Siberian Division of the RAS\\
Novosibirsk, 630090, RUSSIA\\
E-mail: gutman@math.nsc.ru
\par}

\vspace*{5mm}

\noindent
{\it Sem\"en S.~Kutateladze}
\par\smallskip
{\leftskip\parindent\small
\noindent
Sobolev Institute of Mathematics\\
Siberian Division of the RAS\\
Novosibirsk, 630090, RUSSIA\\
E-mail: sskut@math.nsc.ru
\par}

\vspace*{10mm}

\end{document}